\documentclass[10pt]{amsart}
\usepackage{latexsym}
\usepackage{amsmath}
\usepackage{amssymb}
\usepackage{mathrsfs}
\usepackage{graphicx}
\usepackage{xcolor}
\usepackage{bbm}
\usepackage{bbold}
\input epsf
\input xy
\xyoption{all}
\input diagxy

\def\Orth{\mathop{\fam0 Orth}\nolimits}
\def\PSub{\mathop{\fam0 PSub}\nolimits}
\def\Sub{\mathop{\fam0 Sub}\nolimits}
\def\dom{\mathop{\fam0 dom}\nolimits}
\def\On{\mathop{\fam0 On}\nolimits}

\def\Re{\mathop{\fam0 Re}}

\begin{document}

\title[]
{Simultaneous Linear Inequalities:\\ Yesterday and Today}

\author{S.~S. Kutateladze}
\address[]{
Sobolev Institute of Mathematics\newline
\indent 4 Koptyug Avenue\newline
\indent Novosibirsk, 630090\newline
\indent Russia}
\email{
sskut@member.ams.org
}
\begin{abstract}
This is an short overview of  the recent tendencies
in the theory of linear  inequalities that are evoked by Boolean valued
analysis.
\end{abstract}

\date{July 10, 2010}
\thanks{This article bases on a talk
at the opening session of the International Conference
``Order Analysis and Related Problems of Mathematical Modeling,''
Vladikavkaz, July 19--24, 2010, dedicated to the 10th anniversary
of the Vladikavkaz Scientific Center of the Russian Academy of Sciences.}
\maketitle

\section{Agenda}

Linear inequality implies linearity  and order.
When combined, the two produce  an ordered  vector space.
Each linear inequality in the simplest environment  of the sort
is some half-space.
Simultaneity implies many instances and so leads to the
intersections of half-spaces.  These yield polyhedra as well as
arbitrary convex sets, identifying the theory of linear inequalities with
convexity.

Convexity 
reigns in the federation of geometry, optimization, and functional analysis.
Convexity feeds generation, separation, calculus,
and approximation.  Generation appears as duality; separation, as optimality;
calculus,  as representation; and approximation, as stability (cp.~\cite{Harpedonaptae}).

This  talk addresses the origin and the state of the art of the relevant
areas  with a particular emphasis on the Farkas Lemma (cp.~\cite{Farkas_1898}).
Our aim  is to demonstrate how Boolean valued analysis may
be applied to simultaneous linear inequalities with operators.

This particular theme is another illustration of the deep and powerful
technique of ``stratified validity'' which is characteristic of
Boolean valued analysis.

\section{Founding Fathers}

Linearity, inequality, and convexity   stem from the remote
ages (cp.~\cite{History}--
\cite{Encyc}).
However, as the acclaimed pioneers who propounded these ideas and
anticipated their significance for the future
we must rank  the three  polymaths,
{\it Joseph-Louis Lagrange} (January 25, 1736--April 10, 1813),
 {\it Jean Baptiste Joseph Fourier} (March 21, 1768--May 16, 1830), and
{\it Hermann Minkowski} (June 22, 1864--January 12, 1909).

The translator of  the famous elementary textbook of
Lagrange~\cite{Lagrange} Thomas McCormack remarked:

\begin{itemize}
\item[]{\small\it\quad\quad
In both research and exposition, he totally
reversed the methods of his predecessors. They had proceeded in
their exposition from special cases by a species of induction; his
eye was always directed to the highest and most general points of
view; and it was by his suppression of details and neglect of minor,
unimportant considerations that he swept the whole field of analysis
with a generality of insight and power never excelled, adding
to his originality and profundity a conciseness, elegance, and
lucidity which have made him the model of mathematical writers.}
\end{itemize}

The pivotal figure  was Fourier.
Jean-Pierre Kahane wrote in~\cite[pp.~83--84]{Kahane}:

\begin{itemize}
\item[]{\small\it\quad\quad
 He  himself was neglected
for his work on inequalities, what he called ``Analyse ind\'etermi\-n\'ee.''
Darboux considered that he gave the subject an exaggerated importance and
did not publish the papers on this question in his edition of the scientific works
of Fourier. Had they been published, linear programming and convex analysis
would be included in the heritage of Fourier.}
\end{itemize}

David Hilbert lamented the untimely death of Minkowski as follows:\footnote{Cp.~\cite{Hilbert}. Also see~\cite{Corry}.}

\begin{itemize}
\item[]{\small\it\quad\quad
 Since my student years Minkowski was my best, most dependable
friend who supported me with all the depth and loyalty that was so
characteristic of him. Our science, which we loved above all else,
brought us together; it seemed to us a garden full of flowers. In it,
we enjoyed looking for hidden pathways and discovered many a new
perspective that appealed to our sense of beauty, and when one of us
showed it to the other and we marvelled over it together, our joy was
complete. He was for me a rare gift from heaven and I must be grateful
to have possessed that gift for so long. Now death has suddenly torn
him from our midst. However, what death cannot take away is his noble
image in our hearts and the knowledge that his spirit in us continue
to be active.}
\end{itemize}

\section{Environment}

Assume that $X$ is a~real vector space,
$Y$ is a~{\it Kantorovich space\/} also known as a~complete vector lattice
or a Dedekind complete Riesz space.
Let $\mathbb  B:=\mathbb B(Y)$ be the {\it base\/} of~$Y$, i.e., the complete Boolean algebras of
positive projections in~$Y$; and let $m(Y)$ be the  universal completion of~$Y$.
Denote by $L (X,Y)$ the space of linear operators from $X$ to~$Y$.
In case $X$ is furnished with some $Y$-seminorm on~$X$, by $L^{(m)}(X,Y)$ we mean
the {\it space of dominated operators\/} from $X$ to~$Y$.
As usual,  $\{T\le0\}:=\{x\in X \mid Tx\le0\}$; $\ker(T)=T^{-1}(0)$ for $T:X\to Y$.
Also,
$P\in \Sub(X,Y)$ means that $P$ is {\it sublinear}, while
$P\in\PSub(X,Y)$ means that $P$ is {\it polyhedral}, i.e., finitely generated.
The superscript ${}^{(m)}$ suggests domination.

\section{Kantorovich's Theorem}

 Find $\mathfrak X$ satisfying

\[
\xymatrix{
  X\ar[dr]_{B} \ar[r]^{A}
                & W  \ar@{.>}[d]^{\mathfrak X}  \\
                & Y           }
\]

{(\bf 1):}
$(\exists \mathfrak X)\ {\mathfrak X}A=B \leftrightarrow {\ker(A)\subset\ker(B)}.
$

{(\bf 2):} 
{\sl If $W$ is ordered by $W_+$ and $A(X)-W_+=W_+ - A(X)=W$, then}\footnote{Cp.~\cite[p.~51]{SubDif}.}
$$
(\exists \mathfrak X\ge 0)\ {\mathfrak X}A=B \leftrightarrow \{A\le0\}\subset\{B\le 0\}.
$$

\section{The Alternative}

{\sl
Let $X$ be a~$Y$-seminormed real vector space,
with $Y$ a~Kantorovich space.
Assume  that $A_1,\dots,A_N$ and $B$ belong to~$ L^{(m)}(X,Y)$.

Then one and  only one of the following holds:

{\rm(1)} There are   $x\in X$ and $b, b' \in \mathbb B$  such that
$b'\le b$ and
$$
 b'Bx>0, bA_1 x\le0,\dots, bA_N x\le 0.
$$

{\rm(2)} There are   positive orthomorphisms $\alpha_1,\dots,\alpha_N\in\Orth(m(Y))_+$
such that
$
B=\sum\nolimits_{k=1}^N{\alpha_k A_k}.
$
}

\section{Reals: Hidden Dominance}

{\bf Lemma 1.}{\sl\
Let $X$ be a vector space over some subfield $R$ of the reals~$\mathbb{R}$.
Assume that $f$ and $g$ are  $R$-linear functionals on~$X$;
in symbols, $f, g\in X^{\#}:=L(X,\mathbb R)$.

For the inclusion
$$
\{g\le0\}\supset\{f\leq0\}
$$
to hold it is necessary and sufficient that
there be $\alpha\in\mathbb{R}_+$ satisfying
$g=\alpha f$.}

{\scshape Proof.} {\scshape Sufficiency} is obvious.

{\scshape Necessity:} The case of $f=0$ is trivial. If $f\ne0$ then there is
some $x\in X$ such that $f(x)\in \mathbb R$ and~$f(x)>0$. Denote the image $f(X)$ of~$X$
under~$f$ by~$R_0$. Put $h:=g\circ f^{-1}$, i.e. $h\in R_0^{\#}$ is the only solution
for $h\circ f=g$. By hypothesis, $h$ is a positive $R$-linear functional
on~$R_0$. By the Bigard Theorem\footnote{Cp.~\cite[p.~108]{SubDif}.}
$h$ can be extended to a positive homomorphism
$\bar{h}: \mathbb R\to\mathbb R$, since
$R_0-\mathbb R_+=\mathbb R_+-R_0=\mathbb R$.
Each positive automorphism of~$\mathbb R$ is multiplication by a~positive real.
As the sought $\alpha$ we may take  $\bar{h}(1)$.

The proof of the lemma is complete.

\section{Reals: Explicit Dominance}

{\bf Lemma 2.}{\sl\
Let $X$ be an $\mathbb R$-seminormed vector space over some subfield $R$
of~$\mathbb{R}$.
Assume that $f_1,\dots, f_N$ and $g$ are  bounded $R$-linear functionals on~$X$; in symbols, $f_1,\dots, f_N, g\in X^{*}:=L^{(m)}(X,\mathbb R)$.

For the inclusion
$$
\{g\le0\}\supset\bigcap\limits_{k=1}^{N}\{f_k\leq0\}
$$
to hold it is necessary and sufficient that
there be $\alpha_1,\dots,\alpha_N \in\mathbb{R}_+$ satisfying
$
g=\sum\nolimits_{k=1}^{N}\alpha_k f_k.
$
}

\section{Farkas: Explicit Dominance}

{\bf Theorem~1.}{\sl\
Assume  that $A_1,\dots,A_N$ and $B$ belong to~$ L^{(m)}(X,Y)$.

The following are equivalent:

{\rm(1)} Given   $b\in \mathbb B$, the operator inequality $bBx\le 0$
is a consequence of the simultaneous linear operator inequalities
$bA_1 x\le0,\dots, bA_N x\le 0$,
i.e.,
$$
\{bB\le 0\}\supset\{bA_1\le 0\}\cap\dots\cap\{bA_N\le 0\}.
$$

{\rm(2)} There are  positive  orthomorphisms $\alpha_1,\dots,\alpha_N\in\Orth(m(Y))$
such that
$$
B=\sum\limits_{k=1}^N{\alpha_k A_k};
$$
 i.e.,
$B$ lies in the operator convex
conic hull of~$A_1,\dots,A_N$.}

\section{Boolean Modeling}

Cohen's final solution of the problem of the cardinality of the
continuum within ZFC gave rise to the Boolean valued models
by  Scott,  Solovay, and Vop\v enka.\footnote{Cp.~\cite{IBA}.}
Takeuti coined the term {\it Boolean valued
analysis\/}  for applications of the new models
to  analysis.\footnote{Cp.~\cite{Takeuti}.}

Scott forecasted  in 1969:\footnote{Cp.~\cite{Scott}.}

\begin{itemize}
\item[]{\small\it\quad
 We must ask whether there is any interest
in these nonstandard models aside from the independence proof;
that is, do they have any mathematical interest?
The answer must be yes, but we cannot yet give a really good argument.}
\end{itemize}

In 2009 Scott wrote:\footnote{A letter of April 29, 2009 to S.~S. Kutateladze.}

\begin{itemize}
\item[]{\small\it\quad
At the time, I was disappointed that no one took
up my suggestion.  And then I was very surprised
much later to see the work of Takeuti and his
associates.  I think the point is that people have
to be trained in Functional Analysis in order to
understand these models.
I think this is also
obvious from your book and its references.  Alas,
I had no students or collaborators with this kind
of background, and so I was not able to generate
any progress}.
\end{itemize}

\section{Boolean Valued Universe}

Let
${\mathbb B}$
be a~complete Boolean algebra. Given an ordinal
$\alpha$,
put
$$
V_{\alpha}^{({\mathbb B})}
:=\{x \mid
(\exists \beta\in\alpha)\ x:\dom (x)\rightarrow
{\mathbb B}\ \&\ \dom (x)\subset V_{\beta}^{({\mathbb B})}  \}.
$$

The {\it Boolean valued
universe\/}
${\mathbb V}^{({\mathbb B})}$
is
$$
{\mathbb V}^{({\mathbb B})}:=\bigcup\limits_{\alpha\in\On} V_{\alpha}^{({\mathbb B})},
$$
with $\On$ the class of all ordinals.

The truth
value $[\![\varphi]\!]\in {\mathbb B}$ is assigned to each formula
$\varphi$ of ZFC relativized to ${\mathbb V}^{({\mathbb B})}$.

\section{Descending and Ascending}

Given $\varphi$, a~formula of ZFC, and
$y$, a~member of~${\mathbb V}^{{\mathbb B}}$; put
$A_{\varphi}:=A_{\varphi(\cdot,\ y)}:=\{x\mid\varphi (x,\ y)\}$.

The {\it descent\/}
$A_{\varphi}{\downarrow}$
of a~class
$A_{\varphi}$
is
$$
A_{\varphi}{\downarrow}:=\{t\mid t\in {\mathbb V}^{({\mathbb B})} \ \&\  [\![\varphi  (t,\ y)]\!]=\mathbb 1\}.
$$

If
$t\in A_{\varphi}{\downarrow}$,
then
it is said
that
{\it $t$
satisfies
$\varphi (\cdot,\ y)$
inside
${\mathbb V}^{({\mathbb B})}$}.

The {\it descent\/}
$x{\downarrow}$
of $x\in {\mathbb V}^{({\mathbb B})}$
is defined as
$$
x{\downarrow}:=\{t\mid t\in {\mathbb V}^{({\mathbb B})}\ \&\  [\![t\in x]\!]=\mathbb 1\},
$$
i.e. $x{\downarrow}=A_{\cdot\in x}{\downarrow}$.
The class $x{\downarrow}$ is a~set.

If $x$ is a~nonempty set
inside
${\mathbb V}^{({\mathbb B})}$ then
$$
(\exists z\in x{\downarrow})[\![(\exists t\in x)\ \varphi (t)]\!]
=[\![\varphi(z)]\!].
$$

The {\it ascent\/} functor acts in the opposite direction.

\section{The Reals Within}

There is an object
$\mathscr R$
inside
${\mathbb V}^{({\mathbb B})}$ modeling $\mathbb R$, i.e.,
$$
[\![\mathscr R\ {\text{is the reals}}\,]\!]=\mathbb 1.
$$

Let $\mathscr R{\downarrow}$ be the descent of
 the carrier $|\mathscr R|$ of the algebraic system
$
\mathscr R:=(|\mathscr R|,+,\,\cdot\,,0,1,\le)
$
inside ${\mathbb V}^{({\mathbb B})}$.

Implement the descent of the structures on $|\mathscr R|$
to $\mathscr R{\downarrow}$ as follows:
$$
\gathered
x+y=z\leftrightarrow [\![x+y=z]\!]=\mathbb 1;
\\
xy=z\leftrightarrow [\![xy=z]\!]=\mathbb 1;
\\
x\le y\leftrightarrow [\![x\le y]\!]=\mathbb 1;
\\
\lambda x=y\leftrightarrow [\![\lambda^\wedge x=y]\!]=\mathbb 1
\
(x,y,z \in\mathscr R{\downarrow},\ \lambda\in\mathbb R).
\endgathered
$$

{\bf Gordon Theorem.\footnote{Cp.~\cite[p.~349]{IBA}.}}
{\sl $\mathscr R{\downarrow}$
with the descended structures is a~universally complete
vector lattice with base~$\mathbb B(\mathscr R{\downarrow})$
isomorphic to~${\mathbb B}$}.

\section{Proof of Theorem~1}

{(2)$\to$(1):}
If  $B=\sum\nolimits_{k=1}^N{\alpha_k A_k}$ for some positive
$\alpha_1,\dots,\alpha_N$ in~$\Orth(m(Y))$ while  $bA_k x\le 0$ for $b\in\mathbb B$ and
$x\in X$, then
$$
bBx=b\sum\limits_{k=1}^N \alpha_k A_k x=\sum\limits_{k=1}^N \alpha_k bA_k x\le 0
$$
since orthomorphisms commute and projections are orthomorphisms of~$m(Y)$.


{(1)$\to$(2):}
Consider the separated Boolean valued universe $\mathbb{V}^{(\mathbb{B})}$ over the base~$\mathbb B$
of~$Y$.
By~the Gordon Theorem the ascent $Y{\uparrow}$ of~$Y$ is $\mathscr{R}$,
the  reals inside~$\mathbb{V}^{(\mathbb{B})}$.

Using the canonical embedding, we see that $X^{\scriptscriptstyle\wedge}$
is an~$\mathscr R$-seminormed vector space over the standard name  $\mathbb{R}^{\scriptscriptstyle\wedge}$
of the reals~$\mathbb{R}$.

Moreover, $\mathbb{R}^{\scriptscriptstyle\wedge}$ is a~subfield and
sublattice of~${\mathscr R}=Y{\uparrow}$
inside~$\mathbb{V}^{(\mathbb{B})}$.


Put $f_k:=A_k{\uparrow}$ for all $k:=1,\dots, N$  and $g:=B{\uparrow}$.
Clearly,  all $f_1,\dots, f_N,g$ belong to~$(X^{\scriptscriptstyle\wedge})^*$
 inside~${\mathbb V}^{\mathbb B}$.

Define the finite sequence
$$
f:
\{1,\dots,N\}^{\scriptscriptstyle\wedge}\to (X^{\scriptscriptstyle\wedge})^*
$$
  as the
ascent of~$(f_1,\dots,f_N)$.
In other words, the truth values are as follows:
$$
[\![f_{k^{\scriptscriptstyle\wedge}}(x^{\scriptscriptstyle\wedge})=A_k x]\!]={\mathbb 1},\quad [\![g(x^{\scriptscriptstyle\wedge})=Bx]\!]={\mathbb 1}
$$
for all $x\in X$ and $k:=1,\dots,N$.

Put
$$
b:=[\![A_1x\le 0^{\scriptscriptstyle\wedge}]\!]\wedge\dots\wedge
[\![A_N x\le 0^{\scriptscriptstyle\wedge}]\!].
$$
Then $bA_k x\le 0$ for all $k:=1,\dots,N$ and
$bBx\le 0$ by~(1).

 Therefore,
$$
[\![A_1x\le 0^{\scriptscriptstyle\wedge}]\!]\wedge\dots\wedge
[\![A_N x\le 0^{\scriptscriptstyle\wedge}]\!]\le [\![Bx\le 0^{\scriptscriptstyle\wedge}]\!].
$$

In other words,
$$
[\![(\forall k:=1^{\scriptscriptstyle\wedge},\dots,
N^{\scriptscriptstyle\wedge})f_{k}(x^{\scriptscriptstyle\wedge})
\le0^{\scriptscriptstyle\wedge}]\!]
$$
$$
=
\bigwedge_{k:=1,\dots,N} [\![f_{k^{\scriptscriptstyle\wedge}}(x^{\scriptscriptstyle\wedge})\le0^{\scriptscriptstyle\wedge}]\!]
\le [\![g(x^{\scriptscriptstyle\wedge})\le 0^{\scriptscriptstyle\wedge}]\!].
$$

By  Lemma~2 inside~$\mathbb{V}^{(\mathbb{B})}$  and  the maximum principle
of Boolean valued analysis,  there is
a~finite sequence  $\alpha:\{1^{\scriptscriptstyle\wedge},
\dots,N^{\scriptscriptstyle\wedge}\}
\to \mathscr R_+$ inside ${\mathbb V}^{({\mathbb B})}$
satisfying
$$
[\![(\forall x\in X^{\scriptscriptstyle\wedge})\
g(x)
=\sum_{k=1^{\scriptscriptstyle\wedge}}^{N^{\scriptscriptstyle\wedge}}
\alpha(k)  f_k(x)]\!]={\mathbb1}.
$$

Put
$\alpha_k:=\alpha(k^{\scriptscriptstyle\wedge})
\in {\mathscr R}_+{\downarrow}$
for $k:=1,\dots,N$.
Multiplication by an element in~$\mathscr R{\downarrow}$ is an orthomorphism of~$m(Y)$.
Moreover,
$$
B=\sum\limits_{k=1}^N \alpha_k  A_k,
$$
which completes the proof.

\section{ Counterexample: No Dominance }

Lemma~1, describing the consequences of a single inequality,
does not  restrict the class of functionals under consideration.

The analogous version of the Farkas Lemma simply fails for
two simultaneous inequalities in general.

The inclusion
$\{f=0\}\subset\{g\leq 0\}$
equivalent to the inclusion $\{f=0\}\subset\{g=0\}$ does not imply
that $f$ and~$g$ are proportional in the case of an arbitrary subfield
of~$\mathbb R$. It suffices to look at $\mathbb R$ over the rationals $\mathbb Q$,
take some  discontinuous $\mathbb Q$-linear functional on~$\mathbb Q$  and
the identity automorphism of~$\mathbb Q$.

\section{ Reconstruction: No Dominance }

{\bf Theorem 2.}
{\sl
Take $A$ and  $B$ in $L(X,Y)$.
The following are equivalent:

{\rm(1)} $(\exists \alpha\in \Orth(m(Y)))\ B=\alpha A$;

{\rm(2)} There is a projection $\varkappa\in \mathbb B$ such that
$$
\{\varkappa bB\le 0\}\supset\{\varkappa bA\le 0\};
\quad
 \{\neg\varkappa bB\le 0\}\supset\{\neg\varkappa bA\ge 0\}
$$
for all $b\in\mathbb B$.\footnote{As usual, $\neg \varkappa:={\mathbb 1}-\varkappa$.} }

{\scshape Proof.}
Boolean valued analysis reduces the claim to
the scalar case. Applying Lemma~1 twice and writing  down the truth values,
complete the proof.

\section{Interval Operators}

Let $X$ be a~vector lattice.
An~{\it interval operator\/} $\bf T$ from $X$ to~$Y$
is  an~order interval $[\underline{T}, \overline{T}]$
in  $L^{(r)}(X,Y)$, with $\underline{T}\le \overline{T}$.\footnote{Cp.~\cite{Fiedler}.}
The interval equation $\bf B=\mathfrak X \bf A$ has a~{\it weak interval solution\/}
provided that $(\exists \mathfrak X )(\exists A\in \bf A)(\exists B\in \bf B)\ B=\mathfrak X A$.

Given an interval operator $\bf T$ and $x\in X$, put
 $$
 P_{\bf T}(x)=\overline{T}x_+ -\underline{T}x_-.
 $$
\noindent
Call $\bf T$  {\it adapted\/} in case
$\overline{T}-\underline{T}$ is the sum of finitely many disjoint addends.

Put  $\sim(x):=-x$ for all $x\in X$.

\section{Interval Equations}

{\bf Theorem 3.}
{\sl Let $X$ be a vector lattice, and let $Y$ be a~Kantorovich space.
Assume that
 ${\bf A}_1,\dots, {\bf A}_N$ are adapted interval operators and $\bf B$ is an arbitrary
 interval operator in the space of order bounded operators $L^{(r)}(X,Y)$.

The following are equivalent:

{\rm(1)} The interval equation
$$
{\bf B}=\sum_{k=1}^N \alpha_k {\bf A}_k
$$
has a weak interval solution $\alpha_1,\dots,\alpha_N\in\Orth(Y)_+$.

{\rm(2)} For all  $b\in \mathbb B$ we have
$$
\{b\mathfrak B\ge 0\}\supset\{b{\mathfrak A}_1^{\sim}\le 0\}\cap\dots\cap\{b{\mathfrak A}_N^{\sim}\le 0\},
$$
where ${\mathfrak A}_k^{\sim}:=P_{{\bf A}_k}\circ\sim$ for $k:=1,\dots,N$ and
$\mathfrak B:=P_{\bf B}$.
}

\section{Inhomogeneous Inequalities}

{\bf Theorem 4.}
{\sl Let $X$ be a $Y$-seminormed real vector space, with $Y$ a~Kantorovich space.
Assume given some dominated operators
$A_1,\dots,A_N,  B\in L^{(m)}(X,Y)$ and elements $u_1,\dots, u_N,v\in Y$.
The following are equivalent:

{\rm(1)} For all   $b\in \mathbb B$ the  inhomogeneous operator inequality $bBx\le bv$
is a consequence of the
 consistent simultaneous inhomogeneous operator inequalities
 $bA_1 x\le bu_1,\dots, bA_N x\le bu_N$,
i.e.,
$$
\{bB\le bv\}\supset\{bA_1\le bu_1\}\cap\dots\cap\{bA_N\le bu_N\}.
$$

{\rm(2)} There are positive orthomorphisms $\alpha_1,\dots,\alpha_N\in\Orth(m(Y))$
satisfying
$$
B=\sum\limits_{k=1}^N{\alpha_k A_k};\quad v\ge \sum\limits_{k=1}^N{\alpha_k u_k}.
$$
}

\section{Inhomogeneous Matrix Inequalities}

{\bf Theorem 5.\footnote{Cp.~\cite{Olvi}.}}{\sl\
Let $X$ be a~$Y$-seminormed real vector space, with $Y$  a~Kantorovich space.
Assume that
$A\in L^{(m)}(X,Y^s)$, $B\in L^{(m)}(X,Y^t)$, $u\in Y^s$, and $v\in Y^t$,
where $s$ and~$t$ are some naturals.

The following are equivalent:

{\rm(1)} For all   $b\in \mathbb B$  the inhomogeneous operator inequality
$bBx\le bv$ is a~consequence of the consistent inhomogeneous inequality
$bA x\le bu$, i.e., $\{bB\le bv\}\supset\{bA\le bu\}$.

{\rm(2)} There is some $s\times t$ matrix with entries positive orthomorphisms
of~$m(Y)$ such that
$B={\mathfrak X} A$ and~${\mathfrak X}u\le v$
for the corresponding linear operator $\mathfrak X\in L_+(Y^s, Y^t)$.
}

\section{Complex Scalars}

{\bf Theorem 6.}{\sl\
Let $X$ be a~$Y$-seminormed complex vector space, with $Y$ a~Kantorovich space.
Assume given some $u_1,\dots, u_N,v\in Y$ and dominated operators
$A_1,\dots,A_N,  B\in L^{(m)}(X,Y_{\mathbb C})$ from~$X$ into the complexification
$Y_{\mathbb C}:=Y\otimes iY$ of~$Y$.\footnote{Cp.~[3, p.~338].}
Assume further that the  simultaneous inhomogeneous  inequalities
$\vert A_1 x\vert\le u_1,\dots, \vert A_N x\vert\le u_N$ are consistent.
Then the following are equivalent:

\smallskip
$(1)\
\{b\vert B(\cdot)\vert\le bv\}\supset\{b\vert A_1(\cdot)\vert\le bu_1\}\cap\dots\cap\{b\vert A_N(\cdot)\vert\le bu_N\}
$
for all   $b\in \mathbb B$.

{\rm(2)} There are complex orthomorphisms $c_1,\dots, c_N\in\Orth(m(Y)_{\mathbb C})$
satisfying
$$
B=\sum\limits_{k=1}^N{c_k A_k};\quad v\ge \sum\limits_{k=1}^N{\vert c_k\vert u_k}.
$$
}

\section{Inhomogeneous Sublinear Inequalities}

{\bf Lemma 3.}  {\sl
Let $X$ be a real vector space.
Assume that $p_1,\dots, p_N\in \PSub(X):=\PSub(X,\mathbb R)$
and $p\in\Sub(X)$. Assume further that $v, u_1,\dots,u_N\in \mathbb R$
make consistent the simultaneous sublinear inequalities
$p_k(x)\le u_k$, with $k:=1,\dots, N$.

The following are equivalent:

$(1)$\ $\{p\ge v\}\supset\bigcap\limits_{k=1}^{N}\{p_k\leq u_k\};$

$(2)$\ There are $\alpha_1,\dots,\alpha_N\in {\mathbb R}_+$ satisfying

$$
(\forall x\in X)\ p(x)+\sum\limits_{k=1}^{N}\alpha_k p_k(x)\ge0,\quad
\sum\limits_{k=1}^{N}\alpha_k u_k\le -v.
$$
}

{\scshape Proof.}
(2)$\to$(1): If  $x$ is a solution to the
 simultaneous inhomogeneous inequalities
$p_k(x)\le u_k$ with $k:=1,\dots, N$, then
$$
0\le p(x)+\sum\limits_{k=1}^{N}\alpha_k p_k(x)\le
p(x)+\sum\limits_{k=1}^{N}\alpha_k u_k(x) \le p(x)-v.
$$

(1)$\to$(2):   Given $(x,t)\in X\times\mathbb R$, put
$\bar{p}_k(x,t):=p_k(x)-tu_k$, $\bar{p}(x,t):=p(x)-tv$ and
$\tau(x,t):=-t$.
Clearly, $\tau, \bar{p}_1,\dots,\bar{p}_N\in\PSub(X\times\mathbb R)$ and
$\bar{p}\in\Sub(X\times\mathbb R)$.
Take
$$
(x,t)\in \{\tau\le 0\}\cap\bigcap\limits_{k=1}^N\{\bar{p}_k\le 0\}.
$$
If, moreover, $t>0$; then $u_k\ge p_k(x/t)$ for $k:=1,\dots, N$ and so
$p(x/t)\le v$ by hypothesis. In other words $(x,t)\in\{\bar{p}\le 0\}$.
If $t=0$ then take some solution $\bar{x}$ of the
simultaneous inhomogeneous polyhedral inequalities under study.

Since $x\in K:=\bigcap\nolimits_{k=1}^N\{p_k\le 0\}$; therefore,
$p_k(\bar{x}+ x)\le p(x)+ p_k(x)\le u_k$ for all $k:=1,\dots,N$.
Hence, $p(\bar{x}+ x)\ge v$ by hypothesis. So
the sublinear functional $p$ is bounded below on the cone~$K$.
Consequently, $p$ assumes only positive values on~$K$.
In other words, $(x,0)\in\{\bar{p}\le 0\}$. Thus
$$
 \{\bar{p}\ge 0\}\supset\bigcap\limits_{k=1}^{N}\{\bar{p}_k\leq 0\}
$$
and by Lemma~2.2. of~[1]  there are positive reals
$\alpha_1, \dots, \alpha_N, \beta$ such that
for all $(x,t)\in X\times\mathbb R$ we have
$$
\bar{g}(x)+\beta\tau(x)+\sum\limits_{k=1}^N \alpha_k\bar{p}_k(x)\ge 0.
$$
Clearly, the so-obtained parameters $\alpha_1,\dots, \alpha_N$ are what we sought for.
The proof of the lemma is complete.

{\bf Corollary.}
{\sl Let $X$ be an $\mathbb R$-seminormed complex vector space.
Given are $u_1,\dots, u_N,v\in Y$ and
bounded complex linear functionals
$f_1,\dots,f_N,  f\in X^{*}$.
Assume that consistent are the simultaneous inhomogeneous inequalities
$\vert f_1(x)\vert\le u_1,\dots,\vert f_N(x)\vert\le u_N$.
Then the following are equivalent:

{\rm(1)} The inequality $\vert g(x)\vert\le v$
is a consequence of the simultaneous  inequalities
$\vert f_1(x)\vert\le u_1,\dots,\vert f_N(x)\vert\le u_N$,
i.e.
$$
\{\vert g(\cdot)\vert\le v\}\supset \{\vert f_1(\cdot)\vert\le u_1\}\cup\dots\cup\{\vert f_N(\cdot)\vert\le u_N\};
$$

{\rm(2)} There are  $c_1,\dots, c_N\in\mathbb C$
such that}
$$
g=\sum\limits_{k=1}^N{c_k f_k},\quad v\ge \sum\limits_{k=1}^N{\vert c_k\vert u_k}.
$$

{\scshape Proof.}
(2)$\to$(1): If $x\in \bigcap\nolimits_{k=1}^{N}\{\vert f_k(\cdot)\vert\le u_k\}$
then
$$
\vert g(x)\vert=\big\vert \sum\limits_{k=1}^N {c_k f_k(x)}\big\vert\le
\sum\limits_{k=1}^N \vert c_k f_k(x)\vert\le
 \sum\limits_{k=1}^N \vert c_k\vert u_k\le v.
$$

(1)$\to$(2): Consider the realification $X_{\mathbb R}$  of~$X$
and the sublinear functionals
$p(x):=-\Re g(x)$ and $p_k(x):=\vert f_k(x)\vert$, where $k:=1,\dots,N$ and $x\in X_{\mathbb R}$.
Clearly, Lemma~3 applies and
$$
\{p\ge -v\}\supset\bigcap\limits_{k=1}^{N}\{p_k\leq u_k\}.
$$
Hence, there are positive reals
 $\alpha_1,\dots,\alpha_N$ satisfying
$$
(\forall x\in X_{\mathbb R}) -\Re g(x)+\sum_{k=1}^{N}\alpha_k\vert f_k(x)\vert\ge 0;\quad
\sum_{k=1}^{N} \alpha_k u_k\le v.
$$
By subdifferential calculus there are
complexes
$\theta_k, \vert\theta_k\vert=1, k:=1,\dots,N$,
 such that $g=\sum\nolimits\alpha_k\theta_k f_k$.
Put $c_k:=\alpha_k\theta_k$. Obviously,
$\sum\nolimits_{k=1}^{N} \vert c_k\vert u_k=
\sum\nolimits_{k=1}^{N}\alpha_k\vert \theta_k\vert u_k\le v$.
The proof of the corollary is complete.

{\scshape Remark.}
Theorem 6 (which is Theorem 3.1 in~\cite{Farkas} supplied
with a slightly dubious proof)
is a Boolean valued interpretation of the Corollary.

{\bf Theorem 7.}
{\sl
Let $X$ be a $Y$-seminormed real vector space, with $Y$
a~Kantorovich space. Given are some dominated polyhedral sublinear
operators  $P_1,\dots, P_N\in\PSub^{(m)}(X,Y)$ and a~dominated
sublinear operator $P\in\Sub^{(m)}(X,Y)$. Assume further that
$u_1,\dots, u_N,v\in Y$ make consistent the
simultaneous inhomogeneous inequalities $P_1(x)\le u_1,\dots,
P_N(x)\le u_N, P(x)\ge v$.

The following are equivalent:

{\rm(1)} For all  $b\in \mathbb B$ the inhomogeneous sublinear operator
inequality $bP(x)\ge v$  is a~consequence of
the simultaneous inhomogeneous sublinear operator inequalities
$bP_1(x)\le u_1,\dots, bP_N (x)\le u_N$,
i.e.,
$$
\{bP\ge v\}\supset\{bP_1\le u_1\}\cap\dots\cap\{bP_N\le u_N\};
$$

{\rm(2)} There are positive  $\alpha_1,\dots,\alpha_N\in\Orth(m(Y))$
satisfying}
$$
(\forall x\in X)\ P(x)+\sum\limits_{k=1}^N{\alpha_k P_k(x)}\ge 0,\quad
\sum\limits_{k=1}^{N}\alpha_k u_k\le -v.
$$

\section{Lagrange's Principle}

The finite value of the
constrained problem
$$
P_1(x)\le u_1,\dots,P_N(x)\le u_N,\quad P(x)\to\inf
$$
is the value of the unconstrained problem for
an appropriate Lagrangian without any constraint
qualification other that polyhedrality.

The Slater condition  allows us to  eliminate polyhedrality
as well as considering a unique target space.
This is available in a  practically unrestricted generality.\footnote{Cp.~\cite{SubDif}.}

About the new trends relevant to the Farkas Lemma see \cite{Scowscroft}--\cite{Farkas}.

\section{Freedom and Inequality}

Convexity is the theory of linear inequalities in disguise.

Abstraction is the freedom of generalization. Freedom is the loftiest ideal
and idea of man, but it  is  demanding, limited, and vexing.
So is abstraction. So are its instances in convexity, hence,
in simultaneous inequalities.

We definitely feel truth, but we cannot define truth properly.
That is what Alfred Tarski explained to us in the  1930s.\footnote{Cp.~\cite{Tarski}.}
We pursue truth by way of proof, as wittily phrased by
Saunders  Mac Lane.\footnote{Cp.~\cite{Maclane}.}
Mathematics  becomes logic.

The freedom of set theory empowered us with the Boolean valued models
yielding a~lot of surprising and unforeseen visualizations
of the ingredients of mathematics.  Many promising
opportunities are open to modeling  the powerful habits of
reasoning and verification.

Logic  organizes and orders our ways of thinking,
manumitting us from conservatism in choosing the objects and methods
of  research.  Logic of today is  a~fine instrument and
institution of mathematical freedom.
Logic  liberates mathematics
by model theory.

Model theory evaluates and counts truth and proof.
The chase of truth not only leads us close to the truth we pursue but also enables
us to nearly catch up with many other
instances of truth which we were not aware nor even foresaw at the start
of the rally pursuit.  That is what we have learned from
Boolean valued analysis.

Freedom presumes liberty and equality.
Inequality paves way to freedom.

\bibliographystyle{plain}

\end{document}